\documentclass{amsart}
\usepackage{graphicx}
\theoremstyle{definition}

\theoremstyle{remark}

\numberwithin{equation}{section}

\begin{document}

\title[The measure entropy of  LCA with respect to Markov measure]
{The Measure-Theoretical Entropy of  a Linear Cellular
Automata with respect to a Markov Measure} %
\author{Hasan Ak\i n}
\address{Hasan Ak\i n\\
Department of Mathematics\\
Arts and Science Faculty\\
Harran University, \c{S}anliurfa, 63120, Turkey} \email{{\tt
akinhasan@harran.edu.tr}}
\thanks{}%
\renewcommand{\thefootnote}{}
\footnote{2000 \emph{Mathematics Subject Classification}: Primary
28D15; Secondary 37A15.}
\footnote{\emph{Key words and phrases}:
Cellular automata, measure entropy, Markov measure.}

\footnote{A preliminary version of this paper has been presented
to the $19^{th}$ National Turkish Mathematical Symposium,
partially.}


\date{29.08.2006}
\begin{abstract}
In this paper we study the measure-theoretical entropy of the
one-dimensional linear cellular automata (CA hereafter)
$T_{f[-l,r]}$, generated by local rule $f(x_{-l},\ldots,x_{r})=
\sum\limits_{i=-l}^{r}\lambda_{i}x_{i}(\text{mod}\ m)$, where $l$
and $r$ are positive integers, acting on the space of all doubly
infinite sequences with values in a finite ring $\mathbb{Z}_{m}$,
$m \geq 2$, with respect to a Markov measure. We prove that if the
local rule $f$ is bipermutative, then the measure-theoretical
entropy of linear CA $T_{f[-l,r]}$ with respect to a Markov
measure $\mu_{\pi P}$ is\\ $ h_{\mu_{\pi
P}}(T_{f[-l,r]})=-(l+r)\sum\limits_{i,j=0}^{m-1}p_ip_{ij}\text{log}\
p_{ij}.$
\end{abstract}
\maketitle
\section{Introduction}
\par Cellular automata (CA for short), begun by Ulam and von Neumann,
has been systematically studied by Hedlund from purely
mathematical point of view \cite{H}. Hedlund's paper started
investigation of current problems in symbolic dynamics. The study
of such dynamics called CA from the point of view of the ergodic
theory has received remarkable attention in the last few years
(\cite{A1},  \cite{CFMM}, \cite{MM}), because CA have been widely
investigated in a number of disciplines (e.g., mathematics,
physics, computer sciences, etc.). In \cite{Sh}, Shereshevsky has
defined $n$th iteration of a permutative CA and shown that if the
local rule $f$ is right (left) permutative, then its $n$th
iteration also is right (left) permutative.
\par In \cite{MM}, Mass and Martinez have studied the dynamics of Markov measures by
a particular linear cellular automata (LCA). They have reviewed
some results on the evolution of probability measures under CA
acting on a fullshift. In \cite{A5}, the author has studied the
ergodic properties of CA with respect to Markov measures.

\par Although the LCA theory and the entropy
of this LCA have grown up somewhat independently, there are strong
connections between entropy theory and CA theory. We give an
introduction to LCA theory and then discuss the entropy of this
LCA. For a definition and some properties of one-dimensional LCA
we refer to \cite{FLM, H}. The study of the endomorphisms and the
automorphisms (i.e. continuous shift commuting maps, invertible or
non-invertible) of the full shift and its subshifts was initialed
by Hedlund in \cite{H}.

\par It is well known that there are several notions of entropy
(i.e. measure entropy, topological entropy, directional entropy
etc.) of measure-preserving transformation on probability space in
ergodic theory. It is important to know how these notions are
related with each other. The last decade (see. e.g \cite{A1},
\cite{A3}, \cite{MM}), a lot papers are devoted to this subject.

\par In \cite{A1}, the author has compute the measure-theoretical
entropy with respect to uniform Bernoulli measure for the case
$\lambda_i=1$, for all $i\in \mathbb{Z}_m$. The author proved that
the uniform Bernoulli measure is the maximal measure for these
LCA. He has also posed the question whether the maximal measure is
unique.

\par In \cite{A2}, the author has investigated the measure-theoretical directional entropy of
$\mathbb{Z}^{2}$-actions generated by LCA and the shift map with
respect to uniform Bernoulli measure. Let $m, k\in\Bbb N$ and
$\lambda_{-k},\dots,\lambda_k\in\Bbb Z$. Let $\phi\colon\Bbb
Z_m^{\Bbb Z}\to\Bbb Z_m^\Bbb Z$ be a LCA with $\Bbb Z_m=\{0,\dots
,m-1\}$ and local rule $\phi(x)_i=\sum_{j=-k}^k\lambda_j
x_{i+j}\,\pmod m$. Assume $\phi$ is surjective and consider
$v=(v_1,v_2)\in\Bbb Z_+^2$ with $v_2>v_1$. In \cite{A2}, the
author has shown that the directional entropy of $\phi$ with
respect to the uniform Bernoulli measure is given by the formula
$h_v(\phi)=2v_2\text{log}\ m.$

\par In this paper we compute the measure-theoretical entropy of the
one-dimensional LCA $T_{f[-l,r]}$, generated by a bipermutative
local rule\\
 $f(x_{-l},\ldots,x_{r})=
\sum\limits_{i=-l}^{r}\lambda_{i}x_{i}(\text{mod}\ m)$, acting on
the space of all doubly infinite sequences with values in a finite
ring $\mathbb{Z}_{m}$, $m \geq 2$, $\lambda_{i}\in
\mathbb{Z}_{m}$, with respect to a Markov measure $\mu_{\pi P}$.
We show that if the local rule is bipermutative, then we have
$\mu_{\pi P}$ is $ h_{\mu_{\pi
P}}(T_{f[-l,r]})=-(l+r)\sum\limits_{i,j=0}^{m-1}p_ip_{ij}\text{log}\
p_{ij}.$ Also, we prove that if the Markov measure $\mu_{\pi P}$
is uniform, then we get $ h_{\mu_{\pi
P}}(T_{f[-l,r]})=(l+r)\text{log}\ m$.
\par The organization of the paper is as follows: In section 2 we
give the basic formulation of problem to state our main result. In
section 3 we state our main theorem and prove it. In section 4 we
conclude by pointing some further problems.

\section{Preliminaries}

Let ${\mathbb{Z}}_{m}=\{0,1,\ldots, m-1\}$ $(m\geq 2)$ be a ring
of the integers modulo $m$ and $\mathbb{Z}^{\mathbb{Z}}_{m}$ be
the space of all doubly-infinite sequences $x=(x_n)_{n=-\infty
}^{\infty}\in \mathbb{Z}^{\mathbb{Z}}_{m}$ and $x_n\in
{\mathbb{Z}}_{m}$. A CA can be defined as a homomorphism
$\mathbb{Z}^{\mathbb{Z}}_m$ of with product topology. The shift
$\sigma :\mathbb{Z}^{\mathbb{Z}}_{m}\rightarrow
\mathbb{Z}^{\mathbb{Z}}_{m}$ defined by $(\sigma x)_{i} = x_{i+1}$
is a homeomorphism of compact metric space
$\mathbb{Z}^{\mathbb{Z}}_{m}$.
\par A CA is a continuous map, which commutes with
$\sigma$, $T:\mathbb{Z}^{\mathbb{Z}}_{m}\rightarrow
\mathbb{Z}^{\mathbb{Z}}_{m}$ defined by $(Tx)_i = f(
x_{i+l},\ldots, x_{i+r})$, where $f:
\mathbb{Z}^{r-l+1}_{m}\rightarrow \mathbb{Z}_{m}$ is a given local
rule or map. Favati \emph{et al}. \cite{FLM} have defined a local
rule $f$, they have stated that a local rule $f$ is linear
(additive) if and only if it can be written as
\begin{equation}\label{eq1}
f(x_{l},\ldots, x_{r}) = \overset{r}{\underset{i=l}{\sum }}\lambda
_{i}x_{i}(\text{mod}\ m),
\end{equation}
where at least one between $\lambda_{l}$ and $\lambda _{r}$ is
nonzero. We consider 1-dimensional linear cellular automata (LCA)
$T_{f[l, r]}$ determined by the local rule $f$:
\begin{equation}\label{eq2}
(Tx)=(y_n)_{n=-\infty} ^{\infty}, y_n=f(x_{n+l},\ldots ,
x_{n+r})=\overset{r}{\underset{i=l}{\sum }}\lambda
_{i}x_{n+i}(\text{mod}\ m),
\end{equation}
where $\lambda _{l},\ldots, \lambda _{r}\in \mathbb{Z}_{m}$.\\

\par We are going to use the notation $T_{f[l,\ r]}$ for LCA-map
defined in \eqref{eq2} to emphasize the local rule $f$ and the
numbers $l$ and $r$. If the local rule $f$ is given as Eq. (2.1),
then the finite formal power series (\emph{fps} for brevity) $F$
associated with $f$ is defined as
$F(X)=\overset{r}{\underset{i=l}{\sum }}\lambda _{i}X^{-i}$, where
$\lambda _{l},\ldots, \lambda _{r}\in \mathbf{Z}_{m}$ (see
\cite{A3}, \cite{MM} for details). The technique of \emph{fps} is
well known for the study of these problems. In \cite{A3}, the
author has studied the topological entropy of $n$th
iteration of a linear CA by using the \emph{fps}.\\

\par The notion of permutative CA was first introduced by Hedlund
in \cite{H}. If the linear local rule
$f:\mathbb{Z}_{m}^{r-l+1}\to\mathbb{Z}_{m}$ is given in (2.1),
then it is permutative in the $j$th variable if and only if
$gcd(\lambda_j,m)=1$, where $gcd$ denotes the greatest common
divisor. A local rule $f$ is said to be right (respectively, left)
permutative, if $gcd(\lambda_r,m)=1$ (respectively,
$gcd(\lambda_l,m)=1$). It is said that $f$ is bipermutative if it
is both left and right
permutative.\\

\section{The measure entropy of the one-dimensional LCA}
In this section we study the measure entropy of the LCA defined in
Eq. \eqref{eq2} with respect to uniform Markov measure.\\

\par \textbf{Definition 3.1.} Let $(X, \mathcal{B},\mu,T)$ be a
measure-theoretical dynamical system and $\alpha$ be a partition
of $X$. The partition $\alpha$ is called a generator if
$\overset{\infty}{\underset{k=0}{\bigvee }}T^{-k}\alpha
=\mathcal{B}.$\\

\par \textbf{Definition 3.2.} Let $\alpha$ be a partition
of $X$. The quantity
$$
H_{\mu}(\alpha)=-\overset{}{\underset{A\in
\alpha}{\sum}}\mu(A)log\ \mu(A)
$$
is called the entropy of the partition $\alpha$. Let $\alpha$ be a
partition with finite entropy. Then the quantity
$$
h_{\mu}(\alpha,T)=\overset{}{\underset{n\rightarrow\infty}{\lim}}
H_{\mu}(\overset{n-1}{\underset{i=0}{\bigvee}}T^{-i}\alpha)
$$
is called the entropy of $\alpha$ with respect to $T$. One often
writes $h(\alpha,T)$ instead of $h_{\mu}(\alpha,T)$. The quantity
$$
h_{\mu}(T)=\overset{}{\underset{\alpha}{\sup}}\{h_{\mu}(\alpha,T):
\alpha\ \text{is\ a\ partition\ with}\ H_{\mu}(\alpha)<\infty\}
$$
is called the measure-theoretical entropy of $(X,
\mathcal{B},\mu,T)$, the entropy of $T$ (with respect to $\mu$).\\

\par \textbf{Theorem 3.3.} (\cite{W}, Kolmogorov-Sinai Theorem) Let $(X,
\mathcal{B},\mu,T)$ be a measure-theoretical dynamical system and
$\alpha$ a generator with $H_{\mu}(\alpha)<\infty$. Then
$h_{\mu}(T)=h_{\mu}(\alpha,T)$.\\

\par In order to apply entropy theory to the one-dimensional LCA
over ring $\mathbb{Z}_{m}$ ($m\geq 2$) we must define
$\sigma$-algebra $\mathcal{B}$ and the Makov measure
$\mu:\mathcal{B}\rightarrow [0,1]$. In symbolic dynamical system,
it is well known that this $\sigma$-algebra $\mathcal{B}$ is
generated by thin cylinder sets
$$
C = _{a}[j_{0},j_{1},\ldots, j_{s}]_{s+a}=\{x\in
\mathbb{Z}_{m}^{\mathbb{Z}}:x_a=j_{0},\ldots, x_{a+s}=j_{s}\},
$$
where $j_0,j_1,\ldots, j_s \in \mathbb{Z}_{m}$.\\

\par Recall that a subshift of finite type $\sigma:X\rightarrow X$
defined on a space
\begin{equation}
X=\{x\in \mathbb{Z}^{\mathbb{Z}}_m: M_{(x_{n},x_{n+1})}=1,n\in
\mathbb{Z}\}
\end{equation}
for some $m\times m$ matrix $M$ with entries either zero or unity.
Let $P=(p_{(i,j)})$ denote a $m\times m$ stochastic matrix
($p_{(i,j)}\geq 0$, $\sum^{m-1}_{j=0}p_{(i,j)}=1$) with entries
$p_{(i,j)}=0$ iff $M_{(i,j)}=0$, let
$\pi=\{\pi_{0},\pi_{1},\ldots,\pi_{m-1}\}$ be its left
eigenvector. It is well known that $\pi P=\pi$. A pair $\pi$, $P$
defines a set function $\mu_{\pi P}$ on the cylinders of
$\mathbb{Z}_{m}^{\mathbb{Z}}$.  Recall that the associated Markov
measure is defined as follows:
$$
\mu_{\pi
P}(_0[i_0,\ldots,i_{k}]_k)=\pi_{i_{0}}p_{(i_{0},i_{1})}\ldots
p_{(i_{k-1},i_{k})}.
$$

\par See \cite{DGS}, \cite{MM} and \cite{W} for the properties of
the Markov measure.\\

\par Let $\xi$  be the zero-time partition of
$\mathbb{Z}^{\mathbb{Z}}_{m}$: $\xi= \{_{0}[i]: 0\leq i < m\}$,
where $_{0}[i]=\{x\in \mathbb{Z}^{\mathbb{Z}}_{m}: x_0 = i\}$ is a
cylinder set for all $i$, $0 \leq i < m$. So, we can state the
partition $\xi$ as follows:

\begin{equation}\label{eq4}
 \xi= \{_{0}[0], _{0}[1],\ldots, _{0}[m-1] \}.
\end{equation}
Denote by $\mathcal{A}(\xi)$ sub-$\sigma$-algebra of $\mathcal{B}$
generated by the zero-time partition $\xi$ of
$\mathbb{Z}^{\mathbb{Z}}_{m}$.

\par Let us consider a particular case. Assume that the local rule $f$
is bipermutative, so, we have the following Lemma.\\

\par \textbf{Lemma 3.4.} (\cite{A1}, Lemma) Suppose that
$f(x_{-l},\ldots, x_{r})=\overset{r}
{\underset{i=-l}{\sum}}\lambda _{i}x_{i} (mod\ m)$ is a
bipermutative local rule,  and $\xi$ is a partition of
$\mathbb{Z}_{m}^{\mathbb{Z}}$ given in Eq. (3.2), then the
partition $\xi$ is a generator for one-dimensional LCA generated
by $f$.\\

In order to calculate the measure-theoretical entropy of the
one-dimensional LCA $T_{f[-r, r]}$ with respect to uniform Markov
measure we must prove whether
$T_{f[-r, r]}$ is a measure-preserving transformation.\\

\textbf{Proposition 3.5.} Let $T_{f[-l, r]}$ be an one-dimensional
LCA over $\mathbb{Z}_{m}$,  and $f[-l, r]$ be bipermutative local
rule. Then $T_{f[-l, r]}$ is the uniform Markov measure-preserving
transformation.
\begin{proof} Consider a cylinder set
$$
C = _{a}[j_{0}, j_{1},\ldots, j_{s}]_{s+a}=\{x\in
\mathbb{Z}_{m}^{\mathbb{Z}}:x_a^{(0)}=j_{0},\ldots,
x_{a+s}^{(0)}=j_{s}\}.
$$
Then the first preimage of $C$ under $T_{f[-l, r]}$ is the
follows:
\begin{eqnarray*}
T_{f[-l, r]}^{-1}(C)&=&T_{f[-l, r]}^{-1}(\{x\in
\mathbb{Z}_{m}^{\mathbb{Z}}:x_a^{(0)}=j_{0},\ldots,
x_{a+s}^{(0)}=j_{s}\})\\
&=&\overset{}{\underset{(i_{0}, i_{1},\ldots, i_{2r+s})\in
\mathbb{Z}_{m}^{(2r+s+1)}}{\bigcup}}(_{(a-r)}[i_{0}, i_{1},\ldots,
i_{2r+s}]_{a+s+r}),
\end{eqnarray*}
where
$x_a^{(0)}=\overset{r}{\underset{k=-l}{\sum}}\lambda_kx_{a+k}^{(1)}$
(mod $m$) and
$x_{a+s}^{(0)}=\overset{r}{\underset{k=-l}{\sum}}\lambda_kx_{s+a+k}^{(1)}$
(\text{mod}\ $m$). It is clear that,
$$
_{(a-l)}[i_{0},i_{1},\ldots, i_{l+r+s}]_{a+s+r})=\{x\in
\mathbb{Z}_{m}^{\mathbb{Z}}:x_{a-l}^{(1)}=i_{0},\ldots,
x_{a+s+r}^{(1)}=i_{l+r+s}\}.
$$
Then we have
$$
\label{eq1} \mu_{\pi P}(_{a}[j_{0}, j_{1},\ldots,
j_{s}]_{s+a})=\pi_{i_{0}}p_{(j_{0},j_{1})}\ldots
p_{(j_{s-1},j_{s})}.
$$

\begin{eqnarray*}
\mu_{\pi P}(C)&=&\mu_{\pi P}(\{x\in
\mathbb{Z}_{m}^{\mathbb{Z}}:x_a^{(0)}=j_{0},\ldots,
x_{a+s}^{(0)}=j_{s}\})\\
&=&\mu_{\pi P}(T_{f[-l,\ r]}^{-1}(\{x\in
\mathbb{Z}_{m}^{\mathbb{Z}}:x_a^{(0)}=j_{0},\ldots,
x_{a+s}^{(0)}=j_{s}\})\\
&=&\mu_{\pi P}(\overset{}{\underset{i_{0},\ i_{1},\ldots,
i_{(r-l)+s}\in
\mathbb{Z}_{m}^{((r+l)+s+1)}}{\bigcup}}(_{(a-l)}[i_{0},\
i_{1},\ldots,
 i_{(r-l)+s}]_{a+s+r}))\\
&=&m^{l+r}\mu_{\pi P}(_{(a-l)}[i_{0},\ i_{1},\ldots,
 i_{(r+l)+s}]_{a+s+r})\\
&=&m^{l+r}\pi_{i_{0}}p_{(i_{0},i_{1})}\ldots
p_{(i_{r+s-l-1},i_{r+s-l})}\\
&=&m^{-(s+1)}.
\end{eqnarray*}
\end{proof}
\par \textbf{Theorem 3.6.} Let $\mu_{\pi P}$ be a Markov measure
given by the stochastic matrix $P=(p_{ij})$ and the probability
vector $\pi =(p_i)$. Assume that $l$ and $r$ are positive integers
and $gcd(\lambda_l,m)=1$, $gcd(\lambda_r,m)=1$. Then we have
$$
h_{\mu_{\pi
P}}(T_{f[-l,r]})=-(l+r)\sum\limits_{i,j=0}^{m-1}p_ip_{ij}\text{log}\
p_{ij}.
$$
\begin{proof} Now we can calculate the measure entropy of the one-dimensional LCA
by means of the Kolmogorov-Sinai Theorem (\cite{W}, p. 95),
namely, $h_{\mu_{\pi P}}(T_{f[-r,r]})= h_{\mu_{\pi
P}}(T_{f[-r,r]},\xi)$. Let $\xi$  be the zero-time partition of
$\mathbb{Z}^{\mathbb{Z}}_{m}$: $\xi= \{_{0}[i]: 0\leq i < m\}$,
where $_{0}[i]=\{x\in \mathbb{Z}^{\mathbb{Z}}_{m}: x_0 = i\}$ is a
cylinder set for all $i$, $0 \leq i < m$. So, we can state the
partition $\xi$ as follows:
$$
\xi= \{_{0}[0], _{0}[1],\ldots, _{0}[m-1] \}.
$$
Denote by $\mathcal{A}(\xi)$ sub-$\sigma$-algebra of $\mathcal{B}$
generated by the zero-time partition $\xi$ of
$\mathbb{Z}^{\mathbb{Z}}_{m}$. From the definition of entropy we
have
$$
H_{\mu_{\pi P}}(\xi)= -m \mu_{\pi P}(_{0}[i])\text{log}
\mu(_{0}[i])=\text{log}\ m.
$$ From Theorem 3.3 one has
\begin{eqnarray*}
h_{\mu_{\pi P}}(T_{f[-l,r]})&=&-\underset{n\rightarrow \infty
}{\lim } \frac{1}{n}H_{\mu_{\pi
P}}(\overset{n}{\underset{k=0}{\bigvee
}}T_{f[-l, r]}^{-k}\xi )\\
&=& -\underset{n\rightarrow \infty }{\lim } \frac{1}{n}\sum
p_{i_{-nl}}p_{i_{-nl}i_{-nl+1}}\cdots
p_{i_{nr-1}i_{nr}}\times\\
&&\text{log}\ p_{i_{-nl}}p_{i_{-nl}i_{-nl+1}}\cdots
p_{i_{nr-1}i_{nr}}.
\end{eqnarray*}
It is easy to prove by induction that the sum on the right hand
side is equal to
$$
((n(l+r)+1)\sum\limits_{i,j=0}^{m-1}p_ip_{ij}\text{log}\
p_{ij})+\sum\limits_{i=0}^{m-1}p_i\text{log}\ p_{i}.
$$
This implies the result.
\end{proof}

\par A Markov measure on $\mathbb{Z}_{m}^\mathbb{Z}$ is uniform, if measure of any
one-dimensional cylinder is equal to $\frac{1}{m}$ , where $m$ is
a cardinality of $\mathbb{Z}_{m}$. A doubly stochastic matrix is a
matrix $P$ such that $P$ and $P^{tr}$ (transpose) are both
stochastic. If a matrix $P$ is a doubly stochastic then
corresponding Markov measure is a uniform measure. Cardinality of
$\mathbb{Z}_m$ is equal to $m$, so that any doubly stochastic
matrix $P$ of $m\times m$ size will generate uniform Markov measure.\\

\par \textbf{Corollary  3.7.} Let $\mu_{\pi P}$ be the uniform Markov
measure on $\mathbb{Z}_{m}^{\mathbb{Z}}$ and\\
$f(x_{-r},\ldots
,x_{r})=\sum\limits_{i=-r}^{r}\lambda_{i}x_{i}(\text{mod}\ m)$,
where $f[-r,r]$ is bipermutative. Then measure-theoretic entropy
of the one-dimensional LCA $T_{f[-r,r]}$ with respect to $\mu_{\pi
P}$ is equal to $2rlog\ m$.

\begin{proof} It is clear that the partition $\xi\vee T_{f[-r,r]}^{-1}(\xi)$ is
as the following;
$$
\xi\vee T_{f[-r, r]}^{-1}(\xi)=\{_{-r}[i_{-r},\ldots,
i_{r}]_r:i_{-r},\ldots, i_{r}\in \mathbb{Z}_{m} \}.
$$
Because of the uniform Markov measure we get

\begin{eqnarray*}
H_{\mu_{\pi P}}(\xi\vee T_{f[-r,r]}^{-1})&=& -m^{(2r+1)}\mu_{\pi
P}(_{-r}[i_{-r},\ldots , i_{r}]_{r} )\text{log}\ \mu_{\pi
P}(_{-r}[
i_{-r},\ldots , i_{r}]_{r})\\
&=& -m^{(2r+1)}p_{(i_{-r},i_{-r+1})}\cdots
p_{(i_{r-1},i_{r})}\text{log}\
p_{(i_{-r},i_{-r+1})}\cdots p_{(i_{r-1},i_{r})}\\
&=& -m^{(2r+1)} m^{-(2r+1)}\text{log}\ m^{-(2r+1)} =
(2r+1)\text{log}\ m.
\end{eqnarray*}
If we continue, from Lemma 3.4 we have the following results:
\begin{eqnarray*}
H_{\mu_{\pi P}}(\overset{n}{\underset{k=0}{\bigvee
}}T_{f[-r,r]}^{-k}\xi )&=&-m^{(2nr+1)}\mu_{\pi
P}(_{-nr}[i_{-nr},\ldots ,i_{nr}]_{nr} )\times\\
&&log
\mu_{\pi P}(_{-nr}[i_{-nr},\ldots , i_{nr}]_{nr})\\
&=& -m^{(2nr+1)}m^{-(2nr+1)}logm^{-(2nr+1)}\\
&=&(2nr+1)\text{log}\ m.
\end{eqnarray*}
From Theorem 3.3 and Lemma 3.4 we have
$$ h_{\mu_{\pi P}}(T_{f[-r,\ r]})=\underset{n\rightarrow \infty }{\lim }
\frac{1}{n}H_{\mu_{\pi P}}(\overset{n}{\underset{k=0}{\bigvee
}}T_{f[-r,\ r]}^{-k}\xi )= 2r\text{log}\ m.
$$
\end{proof}

\section{Conclusion}

This paper contains the following results: We have found a
generating partition for the one-dimensional LCA generated by a
bipermutative local rule (Lemma 3.4). We have calculated the
measure-theoretical entropy of the one-dimensional LCA with
respect to any Markov measure (Theorem 3.6). This is the first
step toward arbitrary Markov measure classification of
multi-dimensional CA defined on alphabets of composite
cardinality.  In \cite{A1} the author has compute the
measure-theoretical entropy with respect to uniform Bernoulli
measure for the case $\lambda_i=1$, for all $i\in \mathbb{Z}_m$.
The author proved that the uniform Bernoulli measure is the
maximal measure for these LCA. He also posed the question whether
the maximal measure is unique.
\par Thus, where a question raises:\\
Using the Theorem (\cite{W}, Theorem 7.13. (\emph{ii})), can one
calculate the topological entropy of the LCA $T_{f[-r,
r]}:X\rightarrow X$, where $X$ is defined as in Eq. (3.1)? Also it
is open question whether the uniform Markov measure is maximal
measure for $T_{f[-r, r]}:X\rightarrow X$.


\end{document}